\documentclass[12pt,leqno,a4paper]{article}
\usepackage{mathrsfs}

\usepackage{amsfonts}
\usepackage{amsmath, amssymb, amsthm, amsfonts}

\numberwithin{equation}{section}

\def\txt#1{{\textstyle{#1}}}
\def\hf{{\textstyle{\frac12}}}
\def\a{\alpha}\def\b{\beta}
\def\d{{\,\rm d}}
\def\e{\varepsilon}

\def\G{\Gamma} \def\g{\gamma}
\def\k{\kappa}
\def\s{\sigma}
\def\z{\zeta}

\def\={\;=\;}
\def\le{\leqslant}
\def\ge{\geqslant}
\def\zx{\zeta(\hf+ix)}
\def\zt{\zeta(\hf+it)}

\def\D{\Delta}

\def\K{\Bigl({x\over n}\Bigr)^{1/2}\left(K_1(4\pi\sqrt{xn}\,) +
{\pi\over 2}Y_1(4\pi\sqrt{xn}\,)\right)}

\begin{document}
\baselineskip=17pt
\title{\bf \large SOME APPLICATIONS OF LAPLACE TRANSFORMS IN ANALYTIC NUMBER THEORY}
\author{\bf Aleksandar Ivi\'c$^1$}
\date{{\it\small Dedicated to Professor B. Stankovi\'c on the occasion of his 90th birthday}}
\bigskip

\footnotetext[0]{$^1$Serbian Academy of Science and Arts, Knez Mihailova 35,
11000 Beograd, Serbia\\
e-mail: aleksandar.ivic@rgf.bg.ac.rs, aivic\_2000@yahoo.com}
 \maketitle

{\small {\bf Abstract}.
In this overview paper, presented at the meeting DANS14, Novi Sad, July3-7, 2014,
we give some applications of Laplace transforms to analytic number theory. These include
the classical circle and divisor problem, moments of $|\zt|$, and a discussion of two
functional equations connected to a work of Prof. Bogoljub Stankovi\'c.


}
\bigskip

\section{\bf Introduction}

\subsection{Integral transforms}
Integral transforms play an important r\^ole in
Analytic number theory, the
part of Number theory where problems of a number-theoretic nature
are solved by the use of various methods from Analysis.
The most common integral transforms that are used are:
{\it Mellin transforms} (Robert Hjalmar Mellin, 1854-1933), {\it Laplace transforms}
(Pierre-Simon, marquis de Laplace, 1749-1827) and {\it Fourier transforms}
 (Joseph Fourier, 1768-1830). Crudely speaking, suppose that one has
 an integral transform
 \begin{equation}
 F(s) := \int_{\cal I}f(t)K_1(s,t)\d t,
 \end{equation}
where $\cal I$ is an interval, and $K_1(s,t)$ is a suitable kernel function.
If a problem involving the initial function $f(t)$ can be solved by means
of the transforms $F(s)$, then by the inverse transform
\begin{equation}
f(t) := \int_{\cal J}F(s)K_2(s,t)\d s
 \end{equation}
one can obtain information about $f(t)$ itself.
Here $\cal J$ is a suitable contour in
$\Bbb C$ and $K_2(s,t)$ is another kernel function. Naturally the passage from
(1.1) to (1.2) and back requires knowledge about the kernels, the convergence (existence)
of the integrals that are involved, etc.

\subsection{Mellin transforms}

If $f(x)x^{\s-1}\in L(0,\infty)$ and $f(x)$ is of bounded variation in every finite
$x$-interval, then
 \begin{equation}
F(s) := \int_0^\infty f(x)x^{s-1}\d x\qquad(s = \s+it, \; \s,t\in \Bbb R)
 \end{equation}
is the Mellin transform of $f(x)$. The notation of a complex variable $s = \s+it$
in Number Theory originates
with B. Riemann (1826-1866). If (1.3) holds, then under suitable conditions

 \begin{equation}
\hf[f(x+0) + f(x-0)] = \frac{1}{2\pi i}\lim_{T\to\infty}\int_{\s-iT}^{\s+iT}F(s)x^{-s}\d s.
\end{equation}

The relation (1.4) is the  Mellin inversion formula.
Note that if the function
$f(x)$ is continuous, then the left-hand side is simply $f(x)$. As an example, take
$f(x) = x^{s-1}\;(\Re s >0)$. Then $F(s) = \G(s)$, the familiar gamma-function of Euler.
Hence (1.4)  becomes
 \begin{equation}
e^{-z} = \int_{c-i\infty}^{c+i\infty}\G(s)z^{-s}\d s\qquad(c>0, \Re z>0).
\end{equation}
Henceforth we shall use the notation $\int_{c-i\infty}^{c+i\infty}\cdots$ to denote
\[
\lim_{T\to\infty}\int_{c-iT}^{c+iT}\cdots .
\]
The Mellin transform is, because of the presence of $n^{-s}$, particularly useful in dealing with
Dirichlet series
\[
F(s) = \sum_{n=1}^\infty f(n)n^{-s}\qquad(\Re s > \s_0 >0).
\]
A prototype of such series  is the the Riemann zeta-function ($p$ denotes primes)
 \begin{equation}
\zeta(s) = \sum_{n=1}^\infty n^{-s} = \prod_p {(1-p^{-s})}^{-1}\quad(\Re s >1).
\end{equation}
For the values $s$ such that $\Re s \le 1$ (both the series and the product in (1.6)
clearly diverge for $s=1$) $\z(s)$ is defined by analytic continuation
(see e.g., \cite{IVi1} and \cite{IVi2} for an account on $\z(s)$).
\medskip

\subsection{Laplace transforms}

The Laplace transform of $f(x)$ (under suitable conditions on $f(x)$) is

\medskip

\[
{\cal L}\{f(x)\} \equiv F(s) := \int_{0-}^\infty e^{-sx}f(x)\d x \qquad(\Re s >0).
\]
Then ${\cal L}^{-1}\{F(s)\} = f(x)$ is the inverse Laplace transform.
It is unique if e.g., $f(x)$ is continuous.
The inverse Laplace transform can be represented by a complex inversion integral
(so-called Bromwich's inversion integral). This transform is
\[
f(x) = \frac{1}{2\pi i}\lim_{T\to\infty}\int_{\gamma-iT}^{\gamma+iT}e^{sx}F(s)\d s\qquad(x>0)
\]
and $f(x)=0$ for $x<0$. Here $\gamma$ is a real number such that the contour of integration
lies in the region of convergence of $F(s)$.
The Laplace transforms are practical in view of the fast decay factor
$e^{-sx}$, which usually ensures good convergence. For an account on Laplace transforms,
we refer the reader to G. Doetsch's classic works \cite{Doe} and \cite{Doe1}.

\medskip
\subsection{Fourier transforms}

\medskip
Mellin and Laplace transforms are special cases (by a change of variable)
of Fourier transforms (see A.M. Sedletskii \cite{Sed}). This is

\[
{\hat f}(\a) := \int_{-\infty}^\infty f(x)e^{i\a x}\d x\qquad(\a\in\Bbb R)
\]
if $f(x)\in L_1(-\infty, \infty)$. Some sources define ${\hat f}(\a)$ with
the exponential $e^{-i\a x}$ or $e^{2\pi i\a x}$.

\medskip
The inverse Fourier transform is
\[
f(x) = \frac{1}{2\pi}\int_{-\infty}^\infty {\hat f}(\a)e^{-i\a x}\d \a.
\]
This holds for almost all real $x$ if
\[
f(x)\in L_1(-\infty, \infty),\qquad {\hat f}(x)\in L_1(-\infty, \infty).
\]

\medskip
The purpose of this paper is to give first an overview of applications of Laplace transforms
to some solutions of two functional equations, connected to a work of Prof. B. Stankovi\'c,
whose 90th birthday is celebrated this year. Then we shall move to the application
of Laplace transforms to some classical problems of Analytic Number Theory. These include
the circle problem, the divisor problem, the Vorono{\"\i} summation formula for the
divisor function, and moments of $|\zt|$.

\medskip

\section{Solving two functional equations}
\subsection{Prof. Stankovi\'c's work}

\medskip
Walter K. Hayman (1926--\;) in 1967 \cite{Ha} asked:
Is it true that the equation
\begin{equation}
\int_0^\infty \frac{F(wt)}{F(t)}\d t = \frac{1}{1-w}
\end{equation}
has the unique solution $F(t) = ce^t$ (with positive coefficients in
its power series expansion) among the entire functions?

This is a special case of the functional equation
\begin{equation}
\int_0^\infty \frac{F(wt,w)}{F(t,w)}\d t = G(w).
\end{equation}
Both (2.1) and (2.2) were investigated in 1973 by Prof. B. Stankovi\'c \cite{St}
and the speaker \cite{IVi0} (independently), and some partial solutions were obtained.

\footnotetext[0]{$^1$The seminar on Functional Analysis was founded by Prof. B. Stankovi\'c in 1964
and is held at the Department of Mathematics at the University of  Novi Sad  continuously
for 50 years, which is a record that has hardly been attained anywhere. The seminar
was very influential in the scientific formation of many generations of mathematicians.}

In 1972, during the seminar on Functional Analysis in Novi Sad$^1$, Prof. Stankovi\'c noted that
J. Vincze (Budapest) drew his attention
to this problem, which was originally proposed by him and A. R\'enyi (1919-1970)
(after whom the Mathematical Institute of the Hungarian Academy of Sciences was later named).
The author was a member of this Seminar from 1971-1976, and the work that was
subsequently done represents his first research work.

Prof. B. Stankovi\'c proved
\cite{St} several results concerning this
problem. Two of them are the following

\bigskip

{\bf Theorem 1.} {\it
If the integral in} (2.1) {\it exists for a function $F(t)$ when $w_0 \le w < 1$
or $1 < w \le w_0'$, then there does not exist a function $f(t)$  continuous in
a neighborhood $\cal V$ of $t_0>0, f(t_0) \ne0$, and such that}
\[
F(t) = (t-t_0)^kf(t)\qquad(t\in {\cal V}, k \in {\Bbb R}, k> 1 \vee k <-1).
\]

\bigskip

{\bf Theorem 2.} {\it
Let $e^tf(t)$ be a solution of (1) for $G(w) = 1/(1-w), w' \le w < 1$. If $f(t)$ is a
monotonic function of constant sign, then $f(t)$ is a constant.
}

\vfill\break
\medskip
\subsection{The author's work}

\medskip
Under the substitution $F(t) = e^tf(t,w)$ the equation (2.2) becomes

\[
\int_0^\infty e^{-(1-w)t}\frac{f(wt,w)}{f(t,w)}\d t = G(w),
\]
which gives, in the special case when $G(w) = 1/(1-w)$,
\[
\int_0^\infty e^{-(1-w)t}\left\{\frac{f(wt,w)}{f(t,w)}-1\right\}\d t =0.
\]
Therefore we obtain, with
\[
1-w = s, \;G(1-s) = g(s),\; H(s,t) = \frac{f(t-st,1-s)}{f(t,1-s)},
\]
\[
\int_0^\infty e^{-st}H(s,t)\d t = g(s).\leqno(2.3)
\]
The representation (2.3) allows, by the use of the inverse
Laplace transform, to find some families of solutions.

{\bf Theorem 3.} {\it
A solution of the functional equation
\[
\int_0^\infty \frac{F(wt,w)}{F(t,w)}\d t \;=\; \frac{h(w)}{1-w}
\]
is given by
\[
F(t,w)\; =\; t^c \exp\left(\frac{w^c}{h(w)}t\right)\qquad(0<w<1,\; h(w)>0),
\]
and is also valid for $\;w>1\;$ if $\;h(w) <0\;$, where $\,c\,$ is a constant.}

\medskip

\section{Laplace transform in the circle and divisor problem}
\subsection{Introduction}
The circle and divisor problems represent two classical problems
of Analytic Number Theory. Much work was done on them, and the reader
is referred e.g., to \cite{IVi1}, \cite{HI}.

\medskip
The (Gauss) circle problem is the estimation of
\[
P(x) \;:=\; {\sum_{n\le x}}'r(n) - \pi x + 1,
\]
where $x>0$, ${\sum\limits_{n\le x}}'$ means that the last term
in the sum is to be halved if $x$ is an integer and
$r(n) = \sum\limits_{n=a^2+b^2}1$ denotes the number of representations of
$n\;(\in\Bbb N)$ as a sum of two integer squares.

\medskip
The (Dirichlet) divisor problem is the estimation of
\[
\D(x) \;:=\; {\sum_{n\le x}}'d(n) - x(\log x +2\gamma -1) - {\txt{1\over4}},
\]
where $x>0$,
$d(n) = \sum_{\delta|n}1$ is the number of divisors of $n$, and
 $\gamma = -\G'(1) = 0.5772157\ldots $ is Euler's constant.

\medskip
One of the main problems is to determine the value of the constants
$\a, \b$ such that
$$
\a = \inf\Bigl\{ \; a\; |\; P(x) \ll x^a\;\Bigr\},\quad
\b = \inf\Bigl\{ \; b\; |\; \D(x) \ll x^b\;\Bigr\}.
$$
It is known that
$$
1/4 \le \a \le 131/416 = 0,314903\ldots,\quad 1/4 \le \b \le 131/416 = 0,314903\ldots .
$$
It is also generally conjectured that $\a = \b = 1/4$, but this seems very difficult
to prove.

\subsection{Mean square results}

The author in two papers \cite{IVi3}
 investigated the Laplace transform of $P^2(x)$ and $\D^2(x)$.
The key tools are the explicit formulas of
G.H. Hardy (1916) \cite{Har}:
\begin{equation}
P(x) = x^{1/2}\sum_{n=1}^\infty r(n)n^{-1/2}J_1(2\pi\sqrt{xn}),
\end{equation}
and G.F. Vorono{\"\i} (1904) \cite{Vor}:
\[
\D(x) = -\frac{2\sqrt{x}}{\pi}\sum_{n=1}^\infty \frac{d(n)}{\sqrt{n}}
\left\{K_1(4\pi\sqrt{xn}\,)+ \frac{\pi}{2}Y_1(4\pi\sqrt{xn}\,)\right\}.
\]
Here $J_1, K_1, Y_1$ are the familiar Bessel functions (see e.g., Lebedev's
monograph \cite{Leb} and Watson's classic \cite{Wat}
for definitions and properties of Bessel functions),
and the above series for $P(x)$ and $\D(x)$ are boundedly,
but not absolutely convergent. We have
\[
J_p(z) = \sum_{k=0}^\infty \frac{(-1)^k(z/2)^{p+2k}}{k!\G(p+k+1)}
\qquad(p\in \Bbb R, z\in \Bbb C),
\]
and more complicated power series expansions for $K_1, Y_1$.

All Bessel functions have asymptotic expansions involving sines,
cosines and negative powers of $z$, to any degree of accuracy.

{\bf Remark 1}. Note that the
 case of the Laplace transform of $P(x)$ and $\D(x)$ is much easier.
 Namely, using (3.1) and
 \[
 \int_0^\infty e^{-sx}x^{\nu/2}J_\nu(2\sqrt{ax}\,)\d x =
 e^{-a/s}a^{\nu/2}s^{\nu-1}\qquad(\Re s>0,\,\Re \nu > -1)
 \]
one obtains immediately
\[
\int_0^\infty e^{-sx}P(x)\d x = \pi s^{-2}\sum_{n=1}^\infty r(n)e^{-\pi^2n/s}
\qquad(\Re s >0),
\]
and an analogous formula holds for $\D(x)$. We have the following results.

\bigskip

{\bf Theorem 4}. {\it For any given $\e>0$
\[
\int\limits_0^\infty P^2(x)e^{-x/T}\d x = \frac{1}{4}\left(\frac{T}{\pi}\right)^{3/2}
\sum_{n=1}^\infty r^2(n)n^{-3/2} - T + O_\e(T^{2/3+\e}).
\]
}

\bigskip
{\bf Theorem 5}. {\it For any given $\e>0$
\[
\int\limits_0^\infty \D^2(x)e^{-x/T}\d x = \frac{1}{8}\left(\frac{T}{\pi}\right)^{3/2}
\sum_{n=1}^\infty d^2(n)n^{-3/2} + TP_2(\log T) + O_\e(T^{2/3+\e}).
\]
}

\bigskip
{\bf Remark 2}. The standard notation $f(x) = O_\e(g(x))$
means that $|f(x)| \le Cg(x)$ for $C = C(\e)>0, g(x) >0 $ and $x\ge x_0(\e)$.

\bigskip
{\bf Remark 3}.
In Theorem 5, $P_2(x) = a_0x^2 + a_1x + a_2\;(a_0>0)$ with
effectively computable $a_0, a_1, a_2$. Moreover, the series over $n$
in both formulas are both
absolutely convergent, since $r(n) = O_\e(n^\e), d(n) = O_\e(n^\e)$.

\bigskip
{\bf Remark 4}. Mean square formulas for $P(x)$ and $\D(x)$ are a classical problem
in Analytic Number Theory. The best known results are due to W.G. Nowak \cite{Now}
and Lau--Tsang \cite{Lau}, respectively. They are
\[
\int_0^T P^2(x)\d x \= AT^{3/2} + O(T\log^{3/2}T\log\log T)\qquad(A>0)
\]
and
\[
\int_0^T \D^2(x)\d x \= BT^{3/2} + O(T\log^3T\log\log T)\qquad(B>0).
\]
There is an analogy with the formulas of Theorem 4 and Theorem 5, but due
primarily to the presence of the factor $e^{-x/T}$, the formulas of Theorem 4
and Theorem 5 are much more precise. For an account on $\D(x)$ and the mean square
of $|\zt|$, see K.-M. Tsang \cite{Ts}.

\medskip
The key formula for the proof of Theorem 4 is
the Laplace transform
\begin{align*}
&\int\limits_0^\infty e^{-st}tJ_1(a\sqrt{t}\,)J_1(b\sqrt{t}\,)\d t\\&
= \exp\left(-\frac{a^2+b^2}{4s}\right)(4s^3)^{-1}\left\{2abI_0\bigl(\frac{ab}{2s}\bigr)
- (a^2+b^2)I_1\bigl(\frac{ab}{2s}\bigr)\right\}.
\end{align*}
This is valid for $\Re s >0; a,,b\in\Bbb R$. The asymptotic expansion for
the Bessel function $I_\nu(x)$, for fixed $|x| \ge 1$, is (first three terms only)
\[
I_\nu(x) = \frac{e^x}{\sqrt{2\pi x}}\left\{1 - \frac{4\nu^2-1}{8x} +
\frac{(4\nu^2-1)(4\nu^2-9)}{128x^2} + O\left(\frac{1}{|x|^{3}}\right)\right\}.
\]
Namely when we square (3.1) we are led to the integrals of the above type
with $s= 1/T, a = \sqrt{2\pi}m, b =  \sqrt{2\pi}n; \,m,n\in\Bbb N$.

\medskip
The arithmetic part of the proof of Theorem 4 and Theorem 5
is based on the formulas of
F. Chamizo \cite{Cha} and Y. Motohashi \cite{Mo1}, respectively. These are

\begin{align*}
\sum_{n\le x}r(n)r(n+h) & = \frac{(-1)^h8x}{h}\sum_{d|h}(-1)^d d + O_\e(x^{2/3+\e}),\\
\sum_{n\le x}d(n)d(n+h) & = x\sum_{i=0}^2(\log x)^i\sum_{j=0}^2c_{ij}\sum_{d|h}
\frac{(\log d)^j}{d} + O_\e(x^{2/3+\e}).
\end{align*}

The point is that $h$, called ``the shift parameter'', is not necessarily fixed,
but may vary with $x$.
These formulas hold uniformly for $1 \le h \le x^{1/2}$, where the $c_{ij}$'s
are absolute constants, and $c_{22} = c_{21} = 0, c_{20} = 6\pi^{-2}$. The exponents
$2/3+\e$ in the above formulas
account essentially for the same exponents in Theorem 4 and Theorem 5.

\medskip
\section{The Vorono{\"\i} summation formula via Laplace transforms}

\medskip

The author \cite{IVi21}
proved the classical Vorono{\"\i} formula
\begin{equation}
\D(x) = -\frac{2\sqrt{x}}{\pi}\sum_{n=1}^\infty \frac{d(n)}{\sqrt{n}}
\left\{K_1(4\pi\sqrt{xn}\,)+ \frac{\pi}{2}Y_1(4\pi\sqrt{xn}\,)\right\}
\end{equation}
by the use of Laplace transforms. Although there are several proofs
of this important result in the literature (see e.g., Chapter 3 of \cite{IVi1}),
this was the first proof of it by means of Laplace transforms.
We shall present now a sketch of the proof.
Denote the right-hand side of (4.1) by $f(x)$. It can be shown that
\begin{equation}
{\cal L}[\D(x)] = {\cal L} [f(x)]\qquad(x>0).
\end{equation}

Suppose that $x_0 \not \in \Bbb N$. Then both $\D(x)$ and $f(x)$ are
continuous at $x = x_0$. Hence by the uniqueness theorem for Laplace
transforms  it follows that (4.1) holds
for $x = x_0$. But if $x \in \Bbb N$, then the validity of (4.1) follows
from the validity of (4.2) when $x \not \in \Bbb N$, as shown e.g., by
M. Jutila \cite{Jut}. Thus it suffices to show that (4.2) holds when 
$x \not \in \Bbb N$. To see this note that, for $\Re s>0$,

\begin{align*}
{\cal L}[\D(x)] &= \int_0^\infty\left({\sum_{n\le x}}'d(n)
- x(\log x + 2\gamma - 1) - {1\over4}\right)e^{-sx}\d x\\&
= \sum_{n=1}^\infty d(n)\int_n^\infty e^{-sx}\,\d x + {\log s - \gamma\over
s^2} - {1\over4s}\\&
 = {1\over s}\sum_{n=1}^\infty d(n)e^{-sn} + {\log s - \gamma\over
s^2} - {1\over4s}\\&
 = {1\over2\pi is}\int\limits_{(2)}\zeta^2(w)\G(w)s^{-w}\d w +
{\log s - \gamma\over
s^2} - {1\over4s}\\&
 = {1\over2\pi is}\int\limits_{(1/2)}\zeta^2(w)\G(w)s^{-w}\,\d w
- {1\over4s}.
\end{align*}

 Here we used the well-known Mellin integral
(1.5) and the series representation
\[
\zeta^2(s) = \sum_{k=1}^\infty k^{-s}\sum_{\ell=1}^\infty \ell^{-s}
= \sum_{k\ell=n,n\ge1}^\infty n^{-s}= \sum_{n=1}^\infty d(n)n^{-s}\qquad (\Re s > 1).
\]

Change of summation and integration was justified by
absolute convergence,
and in the last step the residue theorem was used together with
$$
\zeta(s) = {1\over s - 1} + \gamma + \gamma_1(s - 1) + \ldots\;,
\quad \G(s) = 1 - \gamma(s - 1) + \ldots\qquad (s \to 1).
$$
The crucial step is to show that
\begin{align*}
&\K \\&= {1\over4\pi^2ni}\int
\limits_{(1)}\G(w)\G(w - 1)\cos^2\left({\pi w\over2}\right)
\left(2\pi\sqrt{xn\,}\right)^{2-2w}\d w.
\end{align*}
For this we need properties of Bessel functions and
the functional equation for $\zeta(s)$,
proved by  first B. Riemann \cite{Rie} in 1859, namely
\[
\z(s) = \chi(s)\zeta(1 - s),\quad \chi(s) = 2^s\pi^{s-1}\sin\left({\pi s\over2}
\right)\G(1 - s)\quad(s\in\Bbb C),
\]
\[
\chi(s) = \left(\frac{2\pi}{t}\right)^{\s+it-1/2}e^{it+i\pi/4}
\left(1+O\left(\frac{1}{t}\right)\right)
\quad(t \ge t_0>0).
\]
With these ingredients the proof of (4.2) is completed.

\bigskip
\section{Laplace transforms of moments of $|\zt|$}

\subsection{Introduction}
%Moments of $|\zt|$

Let
\begin{equation}
L_k(s) \;:=\; \int_0^\infty |\zx|^{2k}\,e^{-sx}\d x
\qquad(k \in \Bbb N,\, \Re s > 0)
\end{equation}
denote that Laplace transform of $|\zx|^{2k}$. Also let
\begin{equation}
I_k(T) \;:=\; \int_0^T|\zt|^{2k}\d t\qquad(k \in \Bbb N).
\end{equation}
The investigations of the moments $I_k(T)$ is one of the central themes
in the theory of $\z(s)$ (see e.g., the monographs \cite{IVi1} and \cite{IVi2}).
One trivially has
$$
I_k(T) \le e\int_0^\infty|\zt|^{2k}
e^{-t/T}\d t =  eL_k\left({1\over T}\right).
$$
Therefore any nontrivial bound of the form
\begin{equation}
L_k(\s) \;\ll_\e\; \left({1\over\s}\right)^{c_k+\e}\qquad(\s\to 0+,\,
c_k \ge 1)
\end{equation}
gives (with $\s = 1/T$) the bound
\begin{equation}
I_k(T) \;\ll_\e\; T^{c_k+\e}.
\end{equation}
Conversely,  if (5.4) holds, we obtain (5.3) from the identity
$$
L_k\left({1\over T}\right) \= {1\over T}\,\int_0^\infty I_k(t)
{\rm e}^{-t/T}\d t.
$$
For a more detailed discussion on $L_k(s)$ and $I_k(T)$,
see the author's work \cite{IVi5}.

\bigskip
\subsection{The mean square of $|\zt|$}

\medskip
A classical result of H. Kober \cite{Kob} from 1936 says that, as $\s \to 0+$,
$$
L_1(2\s) = {\gamma-\log(4\pi\s)\over2\sin\s} +
 \sum_{n=0}^Nc_n\s^n + O_N(\s^{N+1})
$$
for any given integer $N \ge 1$, where the $c_n$'s are effectively
computable constants.

For complex values of $s$ the function $L_1(s)$
was studied by F.V. Atkinson \cite{Atk} in 1941. More recently
M. Jutila \cite{Jut2} refined Atkinson's method  and proved

{\bf Theorem 6}. {\it One has
$$
L_1(s) = -ie^{{1\over2}is}\left(\log(2\pi)-\gamma + ({\pi\over2}-s)i\right)
+
$$
$$
+ 2\pi e^{-{1\over2}is}\sum_{n=1}^\infty
d(n)\exp(-2\pi ine^{-is}) + \lambda_1(s)
$$
in the strip $0 < \Re s < \pi$, where the function $\lambda_1(s)$ is
holomorphic in the strip $|\Re s| < \pi$. Moreover, in any strip}
$|\Re s| \le \theta$ with $0 < \theta < \pi$, we have
$$
\lambda_1(s) \;= \; O_\e\left((|s|+1)^{-1}\right).
$$

In 1997 M. Jutila \cite{Jut1} gave a discussion on the application of Laplace
transforms to the evaluation of sums of coefficients of certain
Dirichlet series. His work showed how a powerful tool Laplace transforms can
be in Analytic number  theory in a general setting.

\bigskip
{\bf Remark 5.} For $I_1(T)$ one has
\begin{equation}
I_1(T) = \int_0^T|\zt|^2\d t = T\left(\log\frac{T}{2\pi} + 2\gamma-1\right)+ E(T),
\end{equation}
say, where $\g$ is Euler's constant and $E(T)$ is the error term in the
asymptotic formula (5.5). For an account on $E(T)$ see e.g., \cite{IVi1}
and \cite{IVi2}. If
$$
\rho \;:=\;
 \inf\Bigl\{ \; r\; |\; E(x) \ll x^r\;\Bigr\},
 $$
 then it is known (see N. Watt \cite{Watt}) that
 $1/4 \le \rho \le \frac{131}{416} = 0.314903\ldots\;.$
 We note that Kober's formula for $L_1(s)$ is different (and in
 some ways more precise) than the formula (5.5) for $I_1(T)$.

\section{The Laplace transform of  $|\zt|^4$}

\subsection{The explicit formula}

Atkinson \cite{Atk2} obtained the asymptotic formula, as $\s \to 0+$,
\begin{equation}
L_2(\s) = {1\over\s}\left(A\log^4{1\over\s} + B\log^3{1\over\s}
+ C\log^2{1\over\s} + D\log {1\over\s} + E\right) + \lambda_2(\s),
\end{equation}
where
$$
A = {1\over2\pi^2},\,B ={1\over\pi^{2}}
\Bigl(2\log(2\pi) - 6\gamma + {24\zeta'(2)\over\pi^{2}}\Bigr),
$$
and
$$
 \lambda_2(\s) \;\ll_\e\;\left({1\over\s}\right)^{{13\over14}+\e},
$$
and indicated how the exponent 13/14 can be replaced by 8/9.
The author in \cite{IVi4} sharpened this result and proved, 
by means of spectral theory of the
non-Euclidean Laplacian,  the following result.

\bigskip
{\bf Theorem 7.}
{\it Let $0 \le \phi < {\pi\over2}$ be given. Then for
$0 < |s| \le 1$ and $|\arg s| \le \phi$ we have}
\begin{align*}
L_2(s) &= {1\over s}\left(A\log^4{1\over s} + B\log^3{1\over s} +
C\log^2{1\over s} + D\log{1\over s} + E\right)\\&
+\, s^{-{1\over2}}\left\{\sum_{j=1}^\infty \a_j H_j^3(\hf)\Bigl(
s^{-i\k_j}R(\k_j)\G(\hf + i\k_j) + s^{i\k_j}R(-\k_j)\G(\hf -
i\k_j) \Bigr)\right\}
\\&
+ G_2(s).
\end{align*}

\bigskip

{\bf Remark 6.}
We have
$$
R(y) \;:=\; \sqrt{{\pi\over2}}{\Bigl(2^{-iy}{\G({1\over4} -
{i\over2}y) \over\G({1\over4} +
{i\over2}y)}\Bigr)}^3\G(2iy)\cosh(\pi y)
$$
and in the above region $G_2(s)$ is a regular function
satisfying ($C > 0$  is a suitable constant)
$$
G_2(s) \ll \frac{1}{\sqrt{|s|}}\exp\left\{
-{C\log(|s|^{-1}+20)\over(\log\log(|s|^{-1}+20))^{2/3}
(\log\log\log(|s|^{-1}+20))^{1/3}}\right\},
$$
where $f(x) \ll g(x)$ means the same as $f(x) = O(g(x))$.
Here
$$
\Bigl\{\lambda_j = \kappa_j^2 + {1\over4}\Bigr\} \,\cup\, \{0\}\,
\qquad(j = 1,2,\ldots)
$$
denotes the discrete spectrum of the non-Euclidean Laplacian acting
on $\,SL(2,\Bbb Z)\,$--automorphic forms, and
$$
\a_j =
|\rho_j(1)|^2(\cosh\pi\kappa_j)^{-1},
$$
where $\rho_j(1)$ is the
first Fourier coefficient of the Maass wave form corresponding to
the eigenvalue $\lambda_j$ to which the Hecke series $H_j(s)$ is
attached. We note that
$$
\sum_{\kappa_j\le K}\a_jH_j^3(\hf) \ll K^2\log^CK \qquad(C >
0).
$$
See Y. Motohashi \cite{Mot} for a detailed account on the spectral theory
of the non-Euclidean Laplacian and its applications to the moments of $\z(s)$.
\bigskip

{\bf Remark 7.} For $I_2(T)$ (see (5.2)) one has
\begin{equation}
I_2(T) = (a_0\log^4T + a_1\log^3T + a_2\log^2T + a_3\log T + a_4)T + E_2(T).
\end{equation}
We have $a_0 = 1/(2\pi^2)$, proved already by A.E. Ingham \cite{Ing}. For
other coefficients in (6.2) see e.g., the author's paper \cite{IVi20}
and J.B. Conrey \cite{Con}, who independently obtained the coefficients in a
different form.
We also have $A = a_0 = 1/(2\pi^2)$ in Theorem 7, but it is easy to see
that, for $B$ in (6.1), $B\ne a_1$ etc. Theorem 3 of \cite{IVi20} provides the exact
connection between the two sets of coefficients.
 Since the series in Theorem 7 is absolutely convergent,
it is seen that
\begin{equation}
L_2(1/T) = {T}\left(A\log^4{T} + B\log^3{T} +
C\log^2{T} + D\log{T} + E\right) + O(\sqrt{T}).
\end{equation}

\bigskip

{\bf Remark 8.}
For the results on $E_2(T)$ see \cite{IM1} and \cite{IM2}. For example, one
has $E_2(T) \ll T^{2/3}\log^CT$. This is weaker than the $O$-term in (6.2).
\bigskip
% \vfill
% \eject


\begin{thebibliography}{99}


{\footnotesize

\bibitem{Atk}  Atkinson, F.V., The mean value of the zeta-function on
the critical line,  Quart. J. Math. Oxford {\bf 10}(1939), 122-128.

\bibitem {Atk2} F.V. Atkinson, The mean value of the zeta-function on
the critical line, Proc. London Math. Soc. {\bf 47}(1941), 174-200.

\bibitem{Cha} Chamizo, F., Correlated sums of $r(n)$, J. Math. Soc.
Japan  {\bf 51}(1999), 237-252.

\bibitem{Con} Conrey, J.B., A note on the fourth power moment of the Riemann zeta-function,
B. C. Berndt (ed.) et al., in ``Analytic number theory. Vol. 1.
Proceedings of a conference in honor of Heini Halberstam, Urbana, 1995'',
 Birkh\"auser, Prog. Math. {\bf138}(1996), 225-230.

\bibitem{Doe} Doetsch, G., Theorie und Anwendung der Laplace-Transformation,
Grundlehren der Mathematischen Wissenschaften in Einzeldarstellungen, vol. {\bf47}(1937),
 Berlin, Springer, 8+436 pp.

\bibitem{Doe1} Doetsch, G., Handbuch der Laplace-Transformation. Vol. 1.
Theorie der Laplace-Transformation,
Mathematische Reihe, vol. {\bf14}(1950), Basel, Birkhäuser, 15+581 pp.


\bibitem{Ha}  Hayman, W.K., Research problems in
function theory, Athlone Press, University of London, London, 1967.


\bibitem{Har} Hardy, G.H.,  The average order of the arithmetical
functions $P(x)$ and $\D(x)$, Proc. London Math. Soc. (2){\bf15}(1916), 192-213.

\bibitem{HI} Huxley, M.N., Ivi\'c, A., Subconvexity for the Riemann
zeta-function and the divisor problem,
Bulletin CXXXIV de l'Acad\'emie Serbe des Sciences et des
Arts - 2007, Classe des Sciences math\'ematiques et naturelles,
Sciences math\'ematiques No. {\bf32}, pp. 13-32.


\bibitem{Ing} Ingham, A.E.,  Mean-value theorems in the theory of the Riemann zeta-function, Proc. London
Math. Soc. (2) {\bf27}(1926), 273-300.

\bibitem{IVi0}  Ivi\'c, A., On a certain integral equation, Matemati\v cki Vesnik (Belgrade)
    {\bf10(25)}(1973), 259-262.

\bibitem{IVi1}  Ivi\'c, A., The Riemann-zeta function, John Wiley \& Sons,
New York, 1985 (2nd ed. Dover, Mineola, 2003).

\bibitem{IVi2} Ivi\'c, A., The mean values of the Riemann zeta-function, Tata
    Institute of Fundamental Research, Lecture Notes {\bf82},
    Bombay 1991 (distr. Springer Verlag, Berlin etc.).

\bibitem{IVi20} Ivi\'c, A., On the fourth moment of the Riemann
zeta-function, Publs. Inst. Math. (Belgrade) {\bf 57(71)}(1995), 101-110.

\bibitem{IVi21} Ivi\'c, A.,
The Vorono{\"\i} identity via the Laplace transform, The Ramanujan Journal
 {\bf 2}(1998), 39-45.

\bibitem{IVi3} Ivi\'c, A.,   The Laplace transform of the square in
the circle and divisor problems,    Studia Scient. Math. Hungarica
 {\bf32}(1996), 181-205 and ibid. {\bf37}(2001), 391-399.

\bibitem{IVi4} Ivi\'c, A.,  The Laplace transform of the fourth moment of
of the zeta-function,  Univ. Beograd. Publ. Elektrotehn. Fak. Ser. Mat.
{\bf11}(2000), 41-48.

\bibitem{IVi5} Ivi\'c, A., The Laplace and Mellin transforms of powers of the Riemann
zeta-function, International Journal of Mathematics and Analysis
{\bf1}(2), 2006, 113-140.


\bibitem{IM1}  Ivi\'c A. and  Motohashi Y., The mean square of the
error term for the fourth moment of the zeta-function,  Proc. London Math.
Soc. (3){\bf66}(1994), 309-329.


\bibitem {IM2}  Ivi\'c A. and  Motohashi Y.,  The fourth moment of the
Riemann zeta-function,  J. Number Theory  {\bf51}(1995), 16-45.

\bibitem{Jut} Jutila, M., A method in the theory of exponential sums, LNs {\bf80},
Tata Institute of Fundamental Research, Bombay, 1987 (distr. by Springer Verlag,
Berlin etc.)



\bibitem{Jut1} Jutila, M., Mean values of Dirichlet series via Laplace
transforms,
in  ``Analytic Number Theory" (ed. Y. Motohashi), London Math. Soc.
 LNS {\bf247},  Cambridge University Press, Cambridge, 1997, 169-207.

\bibitem{Jut2} Jutila, M., The Mellin transform of the square of Riemann's
zeta-function,  Periodica Math. Hung. {\bf42}(2001), 179-190.

\bibitem{Kob}  Kober, H., Eine Mittelwertformel der Riemannschen Zetafunktion,
 Compositio Math. {\bf3}(1936), 174-189.

\bibitem{Lau} Lau Y. K., Tsang K. M., On the mean square formula of the
error term in the Dirichlet divisor problem. Mathematical
Proceedings of the Cambridge Philosophical Society {\bf146}(2009),  277-287.

\bibitem{Leb} Lebedev, N.N., Special functions and their applications,
Dover Publications, Inc., New York, 1972.

\bibitem{Mo1} Motohashi, Y., The binary additive divisor problem, Ann. Sci. \'Ecole
Normale Sup\'erieure (4){\bf27}(1994), 529-572.

\bibitem{Mot} Motohashi, Y., Spectral theory of the Riemann
zeta-function, Cambridge University Press, Cambridge, 1997.


\bibitem{Now} Nowak, W. G., Lattice points in a circle: an improved mean-square asymptotics,
Acta Arith. {\bf113}(2004), 259-272.


\bibitem{Rie} Riemann, B., \"Uber die Anzahl der Primzahlen unter einer gegebener
Gr\"o\ss e, Monatsber. Akad. Berlin (1859), 671-680.

\bibitem{Sed} Sedletskii, A.M., Fourier Transforms and Approximations,
Gordon and Breach Science Publishers, Amsterdam, 2000, 272 pp.

\bibitem{St}  Stankovi\'c, B., On two integral equations, Mathematical Structures --
Computational Mathematics -- Mathematical Modelling, Papers dedicated to Professor
L. Iliev's 60th birthday, Sofia, 1975, 439-445.





\bibitem{Ts} Tsang, K.-M., Recent progress on the Dirichlet divisor problem
and the mean square of the Riemann zeta-function. Sci. China Math.
{\bf53}(2010),  no. 9, 2561-2572.

\bibitem{Vor} G.F. Vorono{\"\i}, Sur une fonction transcendante et ses applications
\`a la sommation de quelques s\'eries, Ann. \'Ecole Normale (3){\bf21}(1904),
2-7-267 and ibid. 459-533.

\bibitem{Wat} Watson, G.N., A Treatise on the Theory of Bessel Functions,
Second Edition (1995), Cambridge University Press, Cambridge.

\bibitem{Watt}  Watt, N., A note on the mean square of $|\zt|$,
J. London Math. Soc. {\bf82(2)}(2010), 279-294.

}

\end{thebibliography}
\end{document}